\documentclass[a4paper,11pt]{amsart}
\usepackage[latin1]{inputenc}
\usepackage{amsmath,amsfonts,amssymb}
\usepackage{color}
\newtheorem{thm}{Theorem}
\newtheorem{lem}{Lemma}[section]%[theo]{Lemma}
\newtheorem{rem}{Remark}[section]%[theo]{Remark}
\newtheorem{defn}{Definition}[section]%[theo]{Definition}

\newcommand{\C}{\mathbb C}
\newcommand{\R}{\mathbb R}

\everymath{\displaystyle}

\begin{document}

\title[Study of some holomorphic curves]{Study of some holomorphic curves in $\C^3$ and their projection into the complex projectve space $\C P^2$}

\author{
Fathi Haggui and Abdessami Jalled}

\maketitle
\begin{abstract}
We study holomorphic curves $f:\C\longrightarrow \C^3$ avoiding four complex hyperplanes and a real subspace of real dimension four or five in $\C^3$. We show that the projection of $f$ into the complex projective space $\C P^2$ is not necessarily constant.
\end{abstract}

{\large\bf keywords}
Complex projective space, Holomorphic curves, Kobayashi hyperbolicity

\section{Introduction}
 The classical Picard Theorem \cite{GP} (see also \cite{B.D}) states that every holomorphic map from the complex Euclidean space $\C$ to $\C P^1$ that avoids three points, is constant. This Theorem has been extended to higher dimension by M.Green \cite{GR} who provided with examples of complex Kobayashi hyperbolic manifolds. We note that if $H_1,...,H_m$ are complex hyperplanes in $\C P^n$, then they are said to be in general position if $m\geqslant n+1$ and any $(n+1)$ of these hyperplanes are linearly independent. Let us recall the Green Theorem:\\
~~\\
\textbf{Theorem (Green, \cite{GR}).} {Let $C$ be a union of $2n+1$ complex hyperplanes in general position in $\C P^n$. Then, any holomorphic curve $f:\C\rightarrow\C P^n\setminus C$ is constant.}\\
\\
In particular, for $n=2$, any holomorphic curve $f:\C\rightarrow\C P^2\setminus C$ is constant, where $C$ is a union of five complex lines in general position in $\C P^n$.\\ As a direct consequence of the Green Theorem, the canonical projection into the complex projective space $\C P^2$ of any holomorphic map $f:\C\rightarrow\C^3$ which avoids five complex hyperplanes in $\C^3$ is constant, since its image avoids the projections of the five complex hyperplanes, which are complex projective lines in general position in $\C P^2$ (see Lemma \ref{lem}). Our main goal is to study the projection into $\C P^2$ of a holomorphic curve $f:\C\rightarrow\C^3$ which avoids four complex hyperplanes in general position in $\C^3$ and a real subspace $H$ of real dimension four or five and check if the projection remains constant.\\
 Throughout the paper we identify $\R^6$, endowed with its standard complex structure $J_{st}$, to $\C^3.$
 \begin{defn}
  Let $n\geqslant 3$ and let $\mathcal{H}=(H_1,...,H_n)$ be a family of real subspaces of $\R^6$ such that $codim_\R H_j=2$ for $j=1,...,n$. Then $\mathcal{H}$ is said to be in general position if for every 3-tuple $(i,j,k)$ of distinct integers $i,j,k\in\lbrace 1,...,n\rbrace$,
$$Span_\R(H_i^\perp,H_j^\perp,H_k^\perp)=\R^6.$$ 
Here, if $H$ is a real subspace in $\R^6$, then $H^\perp$ denotes the orthogonal complement of $H$ with respect to the Euclidean metric.
\end{defn}
We first study the case of four real dimensional subspaces in $\C^3.$ We have the following
\begin{thm}\label{thm2}~\\
$(i)$ Let $H_1,...,H_n$ be $n$ complex hyperplanes in $\C^3$ in general position ($n\geqslant 5$), then there exists a non constant holomorphic curve $f:\C\rightarrow\C^3$ which avoid these $n$ hyperplanes and $\pi(f)$ is constant.\\
$(ii)$ Let $H_1,H_2,H_3,H_4$ be four complex hyperplanes in $\C^3$. Then there exists a real subspace $H$ of\hspace{2mm}$\R^6$, of real dimension four, such that $(H,H_j,H_k)$ are in general position for all $j\neq k,~~j,k\in\lbrace 1,...,4\rbrace$, and there exists $f:\C\rightarrow\C^3$ holomorphic, such that
$$f(\C)\bigcap\big(\bigcup_{i=1}^4H_i\bigcup H\big)=\emptyset$$
and $\pi(f)$ is non constant.
\begin{rem}\vspace{1mm} Here $\pi$ denotes the canonical projection from $\C^3\setminus\lbrace0\rbrace$ into $\C P^2$ and $\pi(f):=\pi\circ f.$ Notice that $\pi(f)$ is well-defined in Theorem \ref{thm2} (ii) since $f(\C)\subset\C^3\setminus\lbrace0\rbrace.$
In case (i), according to the Green Theorem and to Lemma \ref{lem} (see below), $\pi(f)$ is constant. 
\end{rem}
\end{thm}
We study then the case of a subspace in $\C^3$ of real dimension five. We have the following:
\begin{thm}\label{thm3}
 Let $H_1,H_2,H_3,H_4$ be four complex hyperplanes in $\C^3$ and let $H$ be a real subspace of $\R^6$ of real dimension five. Let $\tilde{H}$ be a complex hyperplane of $\C^3$ such that $\tilde{H}\subset H$. Then:\\
$\left( 1\right)$ If $(\tilde{H},H_j,H_k)$ are in general position for all $j\neq k,~~j,k\in\lbrace 1,...,4\rbrace$, then every holomorphic map $~f:\C\rightarrow\C^3$ such that $f(\C)\bigcap(\bigcup_{i=1}^4H_i\bigcup H)=\emptyset$ is constant.\\
$\left( 2\right)$ If there exist $H_j,H_k$, $j\neq k,~j,k\in\lbrace1,...,4\rbrace,$ such that $(\tilde{H},H_j,H_k)$ are not in general position, then there exists $f:\C\rightarrow\C^3$, holomorphic, such that $f(\C)\bigcap\big(\bigcup_{i=1}^4H_i\bigcup H\big)=\emptyset$ and $\pi(f)$ is non constant.
\end{thm}~\\
\begin{rem}~~\\
$\left( a\right)$ The existence and uniqueness of $\tilde{H}\subset H$ is explained in the proof of Theorem \ref{thm3}.\\
$\left( b\right)$ The condition "$(\tilde{H},H_j,H_k)$ are not in general position" is equivalent to the condition "$dim_\R Span_\R(\tilde{H}^\perp,H_j^\perp,H_k^\perp)=4$".\\
$\left( c\right)$ The fact of considering four complex hyperplanes is an optimal condition (see the end of section two for more details).
\end{rem}
The paper is organized as follows. In the first section, we give some results and properties. In section two, we prove Theorem \ref{thm2}. Finally, in section three, we prove Theorem \ref{thm3}.

\section{Preliminaries and properties}
In 1972, Fujimoto \cite{FU} (see also M.Green\cite{GR} and \cite{SL}) showed a statement that characterizes the image of a holomorphic map  $f:\C\rightarrow\C P^n$ omitting $(n+p)$ hyperplanes in general position. He proved the following\\
\textbf{Theorem (Fujimoto \cite{FU}, Serge Lang \cite{SL} pp 196}).
 Let $f: \C \rightarrow \C P^n$ be holomorphic. Assume that the image of $f$ lies in the complement of $n + p$ hyperplanes in general position, then this image is contained in a complex projective subspace of complex dimension $\leqslant [n/p].$\\
The version of the Green Theorem stated in the introduction is a particular case of the previous Theorem, with $p=n+1.$
~~\\

In the remaining of the paper we will need the following properties satisfied by the canonical projection in $\C P^2$ of a holomorphic curve $f:\C\rightarrow\C^3$. For $H$ a real subspace of $\R^6$, we denote by $H^\star$ the set $H\setminus\lbrace0\rbrace.$ Then, we have the following Lemma 
\begin{lem}\label{lem}
Let $\pi:\C^3\setminus\lbrace0\rbrace\rightarrow\C P^2$ be the canonical projection. Then:
\begin{enumerate}
\item If H is a complex hyperplane in $\C^3$, then $\pi(H^\star)$ is a complex projective line in $\C P^2$.
\item If $f:\C\rightarrow\C^3$ is holomorphic and H is a complex hyperplane in $\C^3$, then
$$f(\C)\cap H=\emptyset\Rightarrow \pi(f)(\C)\cap \pi(H^\star)=\emptyset.$$
\item If $H_1,H_2,H_3$ are complex hyperplanes in general position in $\C^3$, then $\pi(H^\star_1),\pi(H^\star_2),\pi(H^\star_3)$ are in general position in $\C P^2$.
\end{enumerate}
\end{lem}
\textbf{Notation}: if $Z\in\C P^2$, we denote $\left[z_1:z_2:z_3\right]$ its homogeneous coordinates, where \\
$(z_1,z_2,z_3)\in\C^3.$
\begin{proof}~\\
\vspace{2mm}
\textbf{Point (1).} We may assume that $H=\lbrace(z_1,z_2,z_3)\in\C^3/ a_1z_1+a_2z_2+a_3z_3=0\rbrace$, with $a_1,a_2,a_3\in\C,$ $a_3\neq0.$ Then
$$\pi(H^\star)=\lbrace\left[1:z_2:z_3\right]\in\C P^2/a_1+a_2z_2+a_3z_3=0\rbrace \cup\lbrace[0:1:-\frac{a_2}{a_3}]\rbrace $$
$$=\lbrace[1:z:-\frac{a_1+a_2z}{a_3}],~z\in\C\rbrace\cup\lbrace [0:1:-\frac{a_2}{a_3}]\rbrace.$$
We notice that $[0:1:-\frac{a_2}{a_3}]$ corresponds to $[\frac{1}{\infty}:1:-\frac{a_1+a_2\infty}{a_3\infty}].$ Hence $\pi(H^\star)$ is a projective complex line in $\C P^2.$ \\
\vspace{2mm}
\textbf{Point (2).} We first notice that $\pi(f)$ is well defined since, by assumption $f(\C)\cap H=\emptyset$, which implies that $f(\C)\subset\C^3\setminus\lbrace0\rbrace.$ Assume now, to get a contradiction, that $\pi(f)(\C)\cap\pi(H^\star)\neq\emptyset$. Then there are two possibilities.\\
\underline{Case $(\alpha).$} There exists $z\in\C$ and there exists $\lambda\in\C$ such that 
$$\pi(f)(z)=\left[1:\lambda:-\frac{a_1+a_2\lambda}{a_3}\right].$$
Then, there exists $c_z\in\C^*$ such that $f(z)=\big(c_z,\lambda c_z,-\frac{a_1+a_2\lambda}{a_3}c_z\big).$ In particular $a_1f_1(z)+a_2f_2(z)+a_3f_3(z)=0$, where $f=(f_1,f_2,f_3).$ Hence, $f(z)\in H.$ This is a contradiction.\\
\underline{Case $(\beta).$}  There exists $z\in\C$ such that 
$$\pi(f)(z)=\left[0:1:-\frac{a_2}{a_3}\right].$$
Then, there exists $c_z\in\C^*$ such that $f(z)=\big(0,c_z,-\frac{a_2}{a_3}c_z\big)$ and $a_1f_1(z)+a_2f_2(z)+a_3f_3(z)=0.$ We obtain again that $f(z)\in H:$ this is a contradiction.\\
\textbf{Point (3).}\vspace{2mm} Since $H_1,H_2,H_3$ are complex hyperplanes in $\C^3$, then there is a linear change of coordinates such that the hyperplanes are defined by equations
$$\begin{array}{ccc}
&H_1=&\lbrace (z_1,z_2,z_3)\in\C^3/ z_1=0\rbrace,\\
&H_2=&\lbrace (z_1,z_2,z_3)\in\C^3/ z_2=0\rbrace,\\
&H_3=&\lbrace (z_1,z_2,z_3)\in\C^3/ z_3=0\rbrace.\\
\end{array}
$$
Now by projection into $\C P^2$, we get
$$\begin{array}{ccc}
&\pi(H^\star_1)=&\lbrace [0:1:z];~z\in\C\cup\lbrace\infty\rbrace\rbrace\cup [0:0:1],\\
&\pi(H^\star_2)=&\lbrace [1:0:z];~z\in\C\cup\lbrace\infty\rbrace\rbrace\cup [0:0:1],\\
&\pi(H^\star_3)=&\lbrace [1:z:0];~z\in\C\cup\lbrace\infty\rbrace\rbrace\cup [0:1:0].
\end{array}
$$
Hence $\pi(H^\star_1)\cap\pi(H^\star_2)\cap\pi(H^\star_3)=\emptyset$, meaning that $\pi(H^\star_1),\pi(H^\star_2),\pi(H^\star_3)$ are in general position since there is no triple point.

\end{proof}
\section{Proof of Theorem \ref{thm2}}
To prove theorem \ref{thm2}, we need the following Lemma which characterize the image of a holomorphic map $f:\C\rightarrow\C P^n$ avoiding $2n$ complex hyperplanes in general position. This precises the result of H.Fujimoto \cite{FU}, \cite{SL} pp 196.
\begin{defn}
 Let $H_1,...,H_m$, $m\geqslant 2n$, be hyperplanes of $\C P^n$. We call diagonal, a line passing through the two points $\bigcap_{i\in I}H_i$ and $\bigcap_{j\in J}H_j$, where $card(I)=card(J)=n$ and $I\cap J=\emptyset.$
\end{defn}
\begin{lem}\label{FU2}
 Let $H_1,...,H_{2n}$ be $(2n)$ projective hyperplanes in general position in $\C P^n$. Then there are $\frac{1}{2}C_{2n}^n$ diagonals $\Delta_1,....,\Delta_{\frac{1}{2}C_{2n}^n}$ such that for every holomorphic curve $f:\C\longrightarrow \C P^n\setminus\bigcup_{i=1}^{2n} H_i$, there exists $k_f\in\lbrace 1,...,\frac{1}{2}C_{2n}^n\rbrace$ such that $f(\C)\subset \Delta_{k_f}.$
\end{lem}
\vspace{2mm}
\begin{proof} ~The proof is inspired by the Fujimoto Theorem, \cite{SL} pp 196.\\
Let $f:\C\longrightarrow \C P^n$ be holomorphic, such that $f(\C)\bigcap(\bigcup_{i=1}^{2n}H_i)=\emptyset$.\\
Let $L_1$,..., $L_{2n}$ be linear forms defining the hyperplanes $H_1,...,H_{2n}$, namely $H_k=L_k^{-1}\left(\lbrace0\rbrace\right)$ for $k=1,...,2n$. If $f=\left[ f_1:...:f_{n+1}\right]$, we denote 
$$h_k:=H_k(f),~~~~~~k=1,...,2n.$$ 
Let $I=\lbrace 1,\cdots,2n\rbrace$ be the set of indices and $\sim$ be the equivalence relation defined by $i \sim j$ if $h_i / h_j$ is constant. We take a partition of the set of indices according to $\sim$. First, we know that the complement of a given class $S$ has at most n elements (see \cite{SL} pp 197). Hence $S$ has at least n elements and there are at most two classes.\\
\vspace{2mm}
The case of one class is not possible. In fact, There exists $\alpha_2,...,\alpha_{2n}\in\C$ such that 
$$\textbf{(S)} \  \left\{\begin{array}{ccllll}
h_2&=&\alpha_2 h_1 & \\
h_3&=&\alpha_3 h_1&\\
  &\vdots&\\
  h_{2n}&=&\alpha_{2n} h_1 \\ 
 
\end{array}\right.
$$
Hence $f(\C)\subset(\bigcap_{k=2}^{n+1}H_k)\bigcap H_1=\bigcap_{k=1}^{n+1}H_k=\emptyset$, which is impossible. Hence there are exactly two classes $S_1$ and $S_2$.\\
~We know that each of the two classes $S_1,S_2$ contains $n$ elements. Then there exists a permutation $\sigma:\lbrace 1,...,2n\rbrace\rightarrow\lbrace 1,...,2n\rbrace$ such that
$$S_1=\lbrace \sigma(1),...,\sigma(n)\rbrace,~~ ~~S_2=\lbrace \sigma(n+1),...,\sigma(2n)\rbrace.$$
Hence There exists $\alpha_2,...,\alpha_{n},\beta_{n+1},...,\beta_{2n-1}\in\C$ such that $h_1,...,h_{2n}$ satisfy the systems:
$$(S_1)  \left\{\begin{array}{ccllll}
h_{\sigma(2)}&=&\alpha_2 h_{\sigma(1)} & \\
h_{\sigma(3)}&=&\alpha_3 h_{\sigma(1)}&\\
  &\vdots&\\
 h_{\sigma(n)}&=&\alpha_n h_{\sigma(1)} \\  
 
\end{array}\right.~~~~
(S_2) \left\{\begin{array}{ccllll}
h_{\sigma(n+1)}&=&\beta_{n+1} h_{\sigma(2n)} & \\
h_{\sigma(n+2)}&=&\beta_{n+2} h_{\sigma(2n)}&\\
  &\vdots&\\
 h_{\sigma(2n-1)}&=&\beta_{2n-1} h_{\sigma(2n)}\\ 
\end{array}\right.
$$
Hence 
$$  \left\{\begin{array}{ccllll}
f(\C)&\subset&(\bigcap_{k=2}^nH_k)\bigcap H_1&=&\bigcap_{k=1}^nH_k \\
f(\C)&\subset&(\bigcap_{k=n+1}^{2n-1}H_k)\cap H_{2n}&=&\bigcap_{k=n+1}^{2n}H_k\\

\end{array}\right.$$
Then $f(\C)\subset \Delta_\sigma$, where $\Delta_\sigma$ is the unique diagonal (line) passing through the two points $ \bigcap_{k=1}^nH_{\sigma(k)}$ and $\bigcap_{k=n+1}^{2n}H_{\sigma(k)}$.\\
Now the two points, and consequently $\Delta_\sigma$, are completely determined by $S_1=\lbrace\sigma(1),...,\sigma(n)\rbrace$ since $S_2$ is automatically fixed once $S_1$ is chosen. Hence $\Delta_\sigma$ is completely determined by a choice of a partition of $\lbrace 1,...,2n\rbrace$ into two subsets, each of them containing $n$ elements. There are exactly $\frac{1}{2}C_{2n}^n$ such partitions. This proves the Lemma.
\end{proof}
We may prove now Theorem \ref{thm2}.\\
We denote by $z=(z_1,z_2,z_3)$ the coordinates in $\C^3$, where $z_j=x_j+iy_j$, $j=1,2,3$. Hence $(x_1,y_1,x_2,y_2,x_3,y_3)$ denote the coordinates in $\R^6.$\\
~~\\
\textbf{Point (i).} Consider first the case $n=5$. By a linear change of coordinates, we take the hyperplanes $H_1,H_2,H_3,H_4$ and $H_5$ in standard form defined by the following equations

$$\begin{array}{lllc}
&H_1=&\lbrace (z_1,z_2,z_3)\in\C^3/ z_1=0\rbrace,\\
&H_2=&\lbrace (z_1,z_2,z_3)\in\C^3/ z_2=0\rbrace,\\
&H_3=&\lbrace (z_1,z_2,z_3)\in\C^3/ z_3=0\rbrace,\\
&H_4=&\lbrace (z_1,z_2,z_3)\in\C^3/ z_1+z_2+z_3=0\rbrace,\\
&H_5=&\lbrace (z_1,z_2,z_3)\in\C^3/ a_1z_1+a_2z_2+a_3z_3=0\rbrace,~~~~a_j\in\R\setminus\lbrace0\rbrace~~~\forall~j=1,2,3.\\
\end{array}
$$
 By hypothesis $f(\C)\cap\big(\bigcup_{i=1}^5H_i\big)=\emptyset.$ Then there exists $h_1,h_2,h_3:\C\rightarrow\C$, holomorphic, such that
$$f=\left( e^{h_1},e^{h_2},e^{h_3}\right).$$
Moreover, since $\pi(f)(\C)$ omits $\pi (H_i)$ for $i=1,...,5$ (see Lemma \ref{lem}) and $\pi\circ f$ is constant by Green (see \cite{GR}), there exists $(\omega_1,\omega_2,\omega_3)\neq(0,0,0)$ such that for all $z\in\C$, 

$$\left[ e^{h_1(z)}:e^{h_2(z)}:e^{h_3(z)}\right]=\left[\omega_1:\omega_2:\omega_3\right].$$
Therefore
$$\left[1:\frac{e^{h_2(z)}}{e^{h_1(z)}}:\frac{e^{h_3(z)}}{e^{h_1(z)}}\right]=\left[1:\frac{\omega_2}{\omega_1}:\frac{\omega_3}{\omega_1}\right]$$
which implies that
$$  \left\{\begin{array}{cllll}
e^{h_2(z)-h_1(z)}&=&\frac{\omega_2}{\omega_1}\\
~~\\
  e^{h_3(z)-h_1(z)}&=&\frac{\omega_3}{\omega_1} \\
  
\end{array}\right.
 \Rightarrow \left\{\begin{array}{cllll}
e^{h_2(z)}&=&\frac{\omega_2}{\omega_1}e^{h_1(z)}\\
~~\\
  e^{h_3(z)}&=&\frac{\omega_3}{\omega_1}e^{h_1(z)} \\
  
\end{array}\right.$$
Hence $f=(e^{h_1},c_2e^{h_1},c_3e^{h_1})$, with $ 1+c_2+c_3\neq0$, and $f$ is not constant.\\
Essentially the same type of argument works in general. Let $H_1,...,H_n,~~~n\geqslant 5$, be $n$ hyperplanes defined by:
 $$H_k:=\left\lbrace  Z\in\C^3/\sum_{i=1}^3\alpha_i^kz_i=0,~~\alpha_i^k\in\C ,~~1\leqslant k\leqslant n\right\rbrace .$$
By hypothesis $f(\C)\bigcap\big(\bigcup_{i=1}^nH_i\big)=\emptyset,$ then in particular $f(\C)\cap\big(\bigcup_{i=1}^5H_i\big)=\emptyset$ and $f=(e^{h},c_2e^{h},c_3e^{h})$ is not constant, where $h$ is holomorphic from $\C$ to $\C$.\\
Hence, in order that $f$ avoids $H_1,...,H_n$, it is sufficient to choose $c_2,c_3\in\C$ such that for every $k=1,...,n$
$$ \alpha_1^k+\alpha_2^kc_2+\alpha_3^k c_3\neq 0.$$
We point out that what preceeds proves more generally that given a countable set of complex hyperplanes in $\C^3$ passing through the origin, there exists $f:\C\rightarrow\C^3$ not constant and avoiding each hyperplane. This proves Point (i).\\
\textbf{Point (ii).}\vspace{1mm} Let $H_1,H_2,H_3$ and $H_4$ be four complex hyperplanes in general position in $\C^3$. We know that there is a linear change of coordinate such that $H_1,H_2,H_3$ and $H_4$ are defined in standard form by :

$$\begin{array}{lllc}
&H_1=&\lbrace (z_1,z_2,z_3)\in\C^3/ z_1=0\rbrace,\\
&H_2=&\lbrace (z_1,z_2,z_3)\in\C^3/ z_2=0\rbrace,\\
&H_3=&\lbrace (z_1,z_2,z_3)\in\C^3/ z_3=0\rbrace,\\
&H_4=&\lbrace (z_1,z_2,z_3)\in\C^3/ z_1+z_2+z_3=0\rbrace,\\
\end{array}
$$
Then 
$$\begin{array}{ccc}
&H_1^\perp=&Span_\R\big[(1,0,0,0,0,0);(0,1,0,0,0,0)\big],\\
&H_2^\perp=&Span_\R\big[(0,0,1,0,0,0);(0,0,0,1,0,0)\big],\\
&H_3^\perp=& Span_\R\big[(0,0,0,0,1,0);(0,0,0,0,0,1)\big],\\
&H_4^\perp=&Span_\R\big[(1,0,1,0,1,0);(0,1,0,1,0,1)\big].\\
\end{array}$$
We pose now 
$$H=\left\{\begin{array}{cllll}
X_1-X_2&=&0\\
~~\\
 X_1-X_3&=&0 \\
  
\end{array}\right.$$
Then $H^\perp=Span_\R\big[(1,0,1,0,0,0);(1,0,0,0,1,0)\big]$, which of course satisfies the condition $Span_\R(H^\perp,H_j^\perp,H_k^\perp)=\R^6$ for all $j\neq k,~~j,k\in\lbrace 1,...,4\rbrace$.\\
Since $f(\C)\bigcap(\bigcup_{i=1}^4 H_i)=\emptyset,$ then there exists holomorphic functions $f_i:\C\rightarrow\C$, $i=1,2,3$ such that 
$$f=(e^{f_1},e^{f_2},e^{f_3}).$$
Then, by Lemma \ref{lem} (2), $g:=\pi(f)$ satisfies $g(\C)\subset\C P^2\setminus\bigcup_{j=1}^4\pi(H^\star_j)$. Hence $g$ has the following form
\begin{equation}
 g=\left[1:e^{g_2}:e^{g_3}\right],
\end{equation}

where $g_2=f_2-f_1$ and $g_3=f_3-f_1$. According to Lemma \ref{FU2} there exists $\frac{1}{2}C^2_{4}=3$ diagonals $\Delta_{12,34},\Delta_{13,24},\Delta_{14,23}$ such that $g=\pi(f(\C))$ is contained in one of these diagonals, where $\Delta_{ij,kl}$ is the diagonal line passing through $\big(\pi(H^\star_i)\bigcap \pi(H^\star_j)\big)$ and $\big(\pi(H^\star_k)\bigcap \pi(H^\star_l)\big)$.\\
 We recall that 
$$ \begin{array}{cllll}
\pi(H^\star_i)&=&\lbrace \left[z_1:z_2:z_3\right]\in\C P^2: z_i=0\rbrace~For~j=1,2,3,\\
\pi(H^\star_4)&=&\lbrace \left[z_1:z_2:z_3\right]\in\C P^2: z_1+z_2+z_3=0\rbrace.\\

\end{array}
$$

Hence $\Delta_{12,34},\Delta_{13,24},\Delta_{14,23}$ are given by

\begin{equation}\label{diagonal}
 \begin{array}{cllll}
\Delta_{12,34}&=&\lbrace \left[z_1:z_2:z_3\right]\in\C P^2: z_1+z_2=0\rbrace,\\
\Delta_{13,24}&=&\lbrace \left[z_1:z_2:z_3\right]\in\C P^2: z_2+z_3=0\rbrace,\\
\Delta_{14,23}&=&\lbrace \left[z_1:z_2:z_3\right]\in\C P^2: z_1+z_3=0\rbrace.
\end{array}
\end{equation}
Suppose that $g(\C)$ is contained in $\Delta_{12,34}$, the cases $g(\C)\subset \Delta_{13,24}$ or $g(\C)\subset \Delta_{14,23}$ being similar. Then\\
 $e^{g_2}+1=0\Rightarrow e^{g_2}=-1\Rightarrow g=\left[1:-1:e^{g_3}\right]$, where $g_3=f_3-f_1$. Hence 
\begin{equation}\label{Forme de f}
 f=(e^{f_1},-e^{f_1},e^{f_3}).
 \end{equation}
On another hand 
$f(\C)\cap H=\emptyset\Leftrightarrow \forall z\in\C$, $\left\{ \begin{array}{cllll}
Re(e^{f_1(z)})&\neq&0\\
or\\
Re(e^{f_1(z)}-e^{f_3(z)})&\neq&0.\\
\end{array}\right.
$~\\
We pose $f_3=2f_1$, then $f=(e^{f_1},-e^{f_1},e^{2f_1})$ avoids $H$. In fact
 $$Re(e^{f_3})=Re(e^{2f_1})=Re(e^{f_1}e^{f_1})=Re(e^{f_1})^2-Im(e^{f_1})^2.$$
Now if $Re(e^{f_1(z)})=0$ for some $z\in\C$, then $Im(e^{f_1(z)})\neq0$ and consequently $Re(e^{f_3}(z))\neq0.$ Hence $f(\C)\cap H=\emptyset.$\\ 
Finally, $\pi(f)=\left[1:-1:e^{f_1} \right]$ is not constant and $f(\C)\bigcap\big(\bigcup_{j=1}^4H_j\bigcup H\big)=\emptyset$. This concludes the proof of Theorem \ref{thm2}.

\section{Proof of Theorem \ref{thm3}}

Let $H$ be a real subspace of $\C^3$ such that $dim_\R H=5$, then $H$ contains a unique complex hyperplane $\tilde{H}$ of $\C^3$. Indeed, there exists $(a_1,b_1,a_2,b_2,a_3,b_3)\in\R^6\setminus\lbrace 0\rbrace$ such that 

$$ \begin{array}{ccl}
H&=&\left\lbrace (x_1,y_1,x_2,y_2,x_3,y_3)\in\R^6/~\sum_{j=1}^3(a_jx_j+b_jy_j)=0\right\rbrace ~~\\
&=&\left\lbrace  z\in\C^3 /~Re\big(\sum_{j=1}^3(a_j-ib_j)z_j\big)=0\right\rbrace .\\

\end{array} $$ 
Hence $\tilde{H}:=\left\lbrace  z\in\C^3 /~\sum_{j=1}^3(a_j-ib_j)z_j=0\right\rbrace $ is a complex hyperplane in $\C^3$, contained in $H$.\\
\textbf{Point (1).}\vspace{1mm} Assume that $(\tilde{H},H_j,H_k)$ are in general position for some $j\neq k$, $j,k\in\lbrace 1,...,4\rbrace$. Since $\tilde{H}\subset H$, where $\tilde{H}$ is a complex hyperplane of $\C^3$, and 
$$f(\C)\bigcap (\bigcup_{i=1}^4H_i\bigcup H)=\emptyset\Rightarrow f(\C)\bigcap (\bigcup_{i=1}^4H_i\cup \tilde{H})=\emptyset,$$ then it follows from Theorem \ref{thm2} (i) that there is $(c_1,c_2)\in(\C^*)^2$ which satisfies $ 1+c_2+c_3\neq0$ and there exists $h:\C\rightarrow\C$ holomorphic such that
 $$f(z)=(e^{h},c_2e^{h},c_3e^{h}).$$
 On another hand $H:=\lbrace (x_1,y_1,...,x_3,y_3)\in\R^6/~\sum_{j=1}^3(a_ix_i+b_iy_i)=0\rbrace.$ By hypothesis $f(\C)\cap H=\emptyset$ then for every $z\in\C$ we have,

$$ \begin{array}{ccc}
&a_1 Re(e^{h(z)})+a_2 Re(c_2 e^{h(z)})+a_3 Re(c_3 e^{h(z)})&\\
&+b_1 Im(e^{h(z)})+b_2 Im(c_2e^{h(z)})+b_3 Im(c_3e^{h(z)})&\neq 0.\\
\end{array} $$     
Thus, for every $z\in\C$
$$ \begin{array}{ccc}
&Re(e^{h(z)})\big[a_1+a_2Re(c_2)+a_3Re(c_3)+b_2Im(c_2)+b_3Im(c_3)\big]&\\
&+Im(e^{h(z)})\big[b_1+b_2Re(c_2)+b_3Re(c_3)-a_1Im(c_2)-a_3Im(c_3)\big]&\neq 0.\\
\end{array} $$    
We denote
$$ \begin{array}{ccc}
&a:=&\big[a_1+a_2Re(c_2)+a_3Re(c_3)+b_2Im(c_2)+b_3Im(c_3)\big]\\~~\\
&b:=&\big[b_1+b_2Re(c_2)+b_3Re(c_3)-a_1Im(c_2)-a_3Im(c_3)\big]\\
\end{array} $$ 
then 
$$f(\C)\cap H=\emptyset\Leftrightarrow e^{h(\C)}\cap \lbrace(x,y)\in\R^2~/~ax+by=0\rbrace=\emptyset.$$       
However $\lbrace(x,y)\in\R^2~/~ax+by=0\rbrace$ is either a real line or $\R^2$, depending on the values of $a$ and $b$. Then by the little Picard Theorem $e^{h}$ is constant because it avoids an infinite number of points. Hence $h$ is constant and $f$ is then constant. We point out that the projection of $f$ into $\C P^2$ is also constant. \\
\textbf{Point (2).}\vspace{1mm} Suppose there exists $j\neq k$, $j,k\in\lbrace1,...,4\rbrace,$ such that\\ $dim_\R Span_\R(\tilde{H}^\perp,H_j^\perp,H_k^\perp)=4$. Then: 
$$\tilde{H}^\perp\subset Span_\R(H_j^\perp,H_k^\perp).$$
In fact for all $i=1,...,4$, $dim_\R H_i^\perp=2$, then $dim_\R Span_\R(H_i^\perp,H_j^\perp)=4$.\\
Suppose $\tilde{H}^\perp\subset Span_\R(H_1^\perp,H_2^\perp)$ then there exists $\alpha_1,\alpha_2\in\C$ such that 
$$\tilde{H}=\lbrace \alpha_1z_1+\alpha_2z_2=0\rbrace.$$
Since $f(\C)\bigcap(\bigcup_{i=1}^4H_i)=\emptyset$, then by \ref{Forme de f} 
$$f=(e^{f_1},-e^{f_1},e^{f_3}).$$
We take $f_1=c,~c\in\C\setminus\lbrace0\rbrace$, such that $Re(\alpha_1e^c-\alpha_2e^c)\neq0$ and $f_3$ not constant. Then 
 $$f=(c,-c,e^{f_3})$$ 
avoids $\bigcup_{i=1}^4H_i\bigcup H$, and $\pi(f)$ is not constant. 
This concludes the proof of Theorem \ref{thm3}.
\qed

By the end of the paper, we show the optimality of considering four complex hyperplanes. Let $H_1,H_2,H_3$ be three complex hyperplanes in $\C^3$, then there exists $H$ a real hyperplane in $\R^6$ and a complex hyperplane $\tilde{H}$ contained in $H$, $\big(H_1,H_2,H_3,\tilde{H}\big)$ are in general position, and there exists $f:\C\rightarrow\C^3$, holomorphic, such that $f(\C)\cap\big(\bigcup_{j=1}^3H_j\bigcup H\big)=\emptyset$ and $\pi\circ f$ is not constant. In fact:\\
We pose $H=\lbrace x_1+x_2+x_3=0\rbrace$ and $\tilde{H}=\lbrace z_1+z_2+z_3=0\rbrace$, which is clearly contained in $H$. Since $f(\C)\bigcap\big(\bigcup_{j=1}^3H_j\bigcup H\big)=\emptyset$, then $f(\C)\bigcap\big(\bigcup_{j=1}^3H_j\bigcup \tilde{H}\big)=\emptyset$ and $f=(e^{f_1},e^{f_2},e^{f_3}).$ Hence 
$$ g:=\pi(f)=\left[1:e^{g_2}:e^{g_3}\right],$$
where $g_2=f_2-f_1$ and $g_3=f_3-f_1$. By lemma \ref{FU2}, $g:=\pi(f)$ is contained in one of diagonals $\Delta_{12,34},\Delta_{13,24},\Delta_{14,23}$ (see \ref{diagonal}). Suppose $\pi(f)(\C)\subset \Delta_{13,24},$ then 
$$\pi(f)=[1,e^{g_2},-e^{g_2}],~~$$
Hence $f=(1,e^{g_2},-e^{g_2})$ avoids $\big(\bigcup_{j=1}^3H_j\bigcup H\big)$ and $\pi(f)$ is not constant. 

\vskip0.5cm 
{\sl E-mail address}:fathi.haggui@gmail.com\\ 
 jalled.abdessami90@gmail.com\\
%~~\\
Institut pr\'eparatoire aux \'etudes d'ing\'enieur de Monastir
Rue Ibn Eljazzar - 5019 Monastir TUNISIE

\end{document}